\documentclass{amsart}

\usepackage[colorlinks, linkcolor=blue,  citecolor=blue, urlcolor=blue]{hyperref}%
\usepackage{url}
\usepackage{verbatim}
\usepackage{alltt}
\usepackage{graphicx}
\graphicspath{{figures/}}

\newtheorem{theorem}{Theorem}[section]

\newtheorem{definition}[theorem]{Definition}
\theoremstyle{remark}
\newtheorem{remark}[theorem]{Remark}
\newtheorem{example}{Example}

\newcommand{\norm}[1]{\Vert#1\Vert}

\newcommand{\abs}[1]{\vert#1\vert}

\newcommand{\grad}{\nabla}
\newcommand{\bq}{\begin{equation}}
\newcommand{\eq}{\end{equation}}
\newcommand{\R}{\mathbb{R}}

\newcommand{\bO}{\mathcal{O}}

\newcommand{\Sk}{{F}}

\usepackage{multirow}
\newcommand{\pad}{\textrm{pad}}

\newcommand{\customwidth}{6in}
\newcommand{\customheight}{4.6in}

\begin{document}

\title[Adaptive Finite Differences for Nonlinear PDEs]
{Adaptive finite difference methods for nonlinear elliptic  and parabolic partial differential equations with free boundaries}

\author{ Adam M Oberman, and Ian Zwiers}

\date{\today}

\begin{abstract}
Monotone finite difference methods provide stable convergent discretizations of a class of degenerate elliptic and parabolic Partial Differential Equations (PDEs).  These methods are best suited to regular rectangular grids, which leads to low accuracy near curved boundaries or singularities of solutions.    In this article we combine monotone finite difference methods with an adaptive grid refinement technique to produce a PDE discretization and solver which is applied to a broad class of equations, in curved or unbounded domains which include free boundaries.    The grid refinement is flexible and adaptive.  The discretization is combined with a fast solution method, which incorporates asynchronous time stepping adapted to the spatial scale.
The framework is validated on linear problems in curved and unbounded domains.  Key applications include the obstacle problem and the one-phase Stefan free boundary problem.
 \end{abstract}
\maketitle


\section{Introduction}
In this article we numerically approximate a class of nonlinear elliptic and parabolic PDEs  using monotone finite difference methods.   Finite difference methods are most easily implemented on regular, rectangular grids.   In this article we combine the monotone finite difference methods with an adaptive  quadtree grid, resulting in significantly improved accuracy near boundaries.   The effectiveness of the method is demonstrated on the Laplace equation on curved and on unbounded domains.   Key applications are the obstacle problem, and the Stefan Free Boundary problems.

 Using the framework of nonlinear elliptic operators, we can combine the partial differential equation with the boundary conditions (or even free boundaries)  into a \emph{single} degenerate elliptic operator.  This allows us to build adaptive discretizations and solvers using a unified framework, and to experiment with different grid adaptation strategies.


Adaptive finite difference methods have been used in a similar context in a variety of problems, but a not in a framework as general as this one.    A review of data structures and implementation of sparse grids for Partial Differential Equations can be found in~\cite{bungartz2004sparse}. 
Many approaches using finite differences methods combine the popular  level set method for tracking the boundary with a representation of the operator inside the boundary.   
A fourth order adaptive method for the heat equation and stefan equation can be found in~\cite{gibou2005fourth}. Adaptive grids for the Stefan problem were used in~\cite{chen2009numerical}.
  Adaptive grid refinement combined with a level set representation of the free boundary  was used for the Poisson-Boltzmann system in~\cite{helgadottir2011poisson,mirzadeh2011second}.

An advantage of the finite difference implementation and the viscosity solution framework is that the conditioning of the solvers does not break down as the equation becomes degenerate.    For example, fast solvers for the degenerate elliptic Monge-Ampere equation have been built, where the Newton's method solver speed is (nearly) independent of the regularity of the solutions~\cite{ObermanFroeseFast,ObermanFroeseFiltered}.  These problems were solved on a uniform grid using wide stencil finite difference schemes, but the later article extended the problem to Optimal Transportation boundary conditions, where the source domain is irregular, and the target domain is convex~\cite{Benamou2014}.  However the anisotropy of the operator requires wide stencils for monotone discretizations, which are more challenging to implement on an adaptive grid.


In order to work with an adaptive grid, we need a refinement criteria.  We take the point of view that the equation itself should provide this criteria.  By writing the entire problem (including boundary conditions) as a \emph{single degenerate elliptic operator} we are able to produce an effective refinement criteria.

The framework can be used for many purposes, including:
\begin{itemize}
\item Artificial Boundary Conditions for problems in an unbounded domain. We use coarse grids in the exterior, and choose to adapt based on either the residual of the boundary conditions or the distance from a reference point in the domain.
\item Grid adaptation for PDEs on curved domains, using grid based discretizations.
\item Obstacle problems or one phase free boundary problems such as the Stefan problem.
\item Nonlinear iterative methods for stationary problems.   On an adaptive grid, the iterations are asynchronous, so there the nonlinear CFL condition is locally determined.\end{itemize}

\subsection{The framework of  degenerate elliptic operators}\label{sec:degelleqns}

We consider the class of degenerate elliptic equations~\cite{CIL}, which include first order equations, such as the eikonal equation, as well as fully nonlinear PDEs, such as the Monge-Ampere equation, and free boundary problems.
Singularities can be present in the solutions to these equation, in particular at locations near the free boundary or where the equation changes types.  For this reason, weak solutions, are needed, which are the viscosity solutions~\cite{CIL}.  The theory of viscosity solutions is by now well-established.   To prove convergence of the schemes, we require that that uniqueness hold for the underlying PDE.  In most 
cases, this is covered by the standard theory.   Classical solutions of the Stefan problem arise only under limited conditions~\cite{daskalopoulos2005all}.   For the one phase Stefan problem, uniqueness of viscosity solutions is established in~\cite{kim2003uniqueness}.

Let  $\Omega$ be a domain in $\R^n,$ $Du$ and $D^2u$ denote the gradient and Hessian of $u$, respectively,  and let
$F(X,p, r,x)$ be a  continuous real valued function  defined on $\mathbb{S}^n\times \R^n\times \R \times \Omega$,
 $\mathbb{S}^n$ being the space of symmetric $n\times n$ matrices.    Write 
\[
F[u](x) \equiv F(D^2u(x),Du(x),u(x), x).
\]  
\begin{definition}
\label{defn:degell}  The operator $F$ is  \emph{degenerate elliptic} if
\[
F(X,p, r,x) \leq F(Y,p,s,x) \text{ whenever } r \leq s \text{ and } Y \leq X,
\]
where $Y \leq X$ means that $Y - X$ is a nonnegative definite symmetric matrix. 
\end{definition}

If the operator $F$ is degenerate elliptic, then we say the Partial Differential Equation on the domain $\Omega$ 
\[
F[u](x) 
= 0, \quad \text{ for $x$ in } 
\] 
(along with, for example,  Dirichlet boundary conditions,
$ u(x) = g(x)$ , or $x$ on $\partial \Omega$)  is as well.
The initial-boundary value problem for the 
\[
u_t(x,t) + F[u](x,t)
\]
is called degenerate parabolic, when the operator $F$ is degenerate elliptic.

\begin{example}[Examples of degenerate elliptic operators]
The obstacle problem, 
\[
\min( -u_{xx}, u -g(x)) = 0
\] is degenerate elliptic. 
The Hamilton-Jacobi equation
\[
u_t - |u_x| = 0,
\] is degenerate parabolic.  
The equation 
\bq\label{GFbc}
c(x) ( - \Delta u(x) + f(x)) + d(x)(u(x) - g(x))
\eq
is degenerate elliptic, provided $c(x), d(x) \ge 0$. 
\end{example}

\subsection{Elliptic finite difference methods}


The class of finite difference methods (or equations) we focus on are called \emph{elliptic}, \cite{ObermanSINUM}.  They are a special class of monotone finite difference schemes which are automatically stable, and arise from a simple construction.  Consistent elliptic schemes, since they are monotone and stable,  converge, according to the theory presented in \cite{BSnum}.


Finite difference equations can be defined on a general unstructured grid, regarded as a weighted, directed graph.
In our case the adaptive finite difference grid has a natural data structure given by the quadtree, which is discussed below.  But to define monotone schemes, we can consider the abstract setting.
The  unstructured grid on the domain $\Omega$;  is a directed graph consisting of a set of points, $x_i \in \Omega, i = 1,\dots N$, each endowed with a list of neighbors, $N(i) = (i_1, \dots, i_d)$.  
A  \emph{grid  function} is a real-valued function defined on the grid, with values $u_{i} = u(x_{i}).$  
The finite difference operator is represented at each grid point by 
\bq\label{finitediff}
\Sk^i[u] \equiv \Sk^i 
\left(u_{i},  \frac{u_{i} - u_{i_1}} { \abs{x_i - x_{i_1}}}, \dots, \frac{u_{i} - u_{i_d}} { \abs{x_i - x_{i_d}}} \right), \quad  i = 1,\dots N,
\eq
where $\Sk^i(x,y_1, \dots, y_d)$ is a specified, usually nonlinear, function of its arguments.   The list of finite differences in the above expression can be regarded as the gradient of the function on the graph.  The notation
\[
\grad u(x_i) = 
\left(\frac{u_{i} - u_{i_1}} { \abs{x_i - x_{i_1}}}, \dots, \frac{u_{i} - u_{i_d}} { \abs{x_i - x_{i_d}}} \right)
\]
was used in~\cite{manfredi2012nonlinear}, so that we can write
\[
\Sk^i[u] = \Sk^i (u_i, \grad u(x_i) ).
\]
This notation emphasizes the fact that a finite difference operator is local: it depends only on the value at the reference points, and the gradient of the function on the graph.  (Second order finite differences come from combinations of first order differences; higher order differences are not needed).
A \emph{solution} is a grid function which satisfies  $\Sk[u] = 0$ (at all grid points).  
A \emph{boundary} point  can be identified as a grid point with no neighbors, so that  Dirichlet boundary conditions can be  imposed  by setting  $\Sk^i[u] = u_i - g(x_i)$.

We now define degenerate elliptic operators.
\begin{definition}\label{defn:degenElliptic}    
The finite difference operator $\Sk$ is \emph{degenerate elliptic} if each component $\Sk^i(x,y_1, \dots, y_d)$ is nondecreasing in each variable.
\end{definition}
We emphasize that the scheme is a nondecreasing function of $u_i$ and the differences $u_i - u_j$ for neighbors $j$ of $i$.

\begin{remark}We now explain the reason for using degenerate elliptic schemes in this context.  
On a uniform grid, the standard discretization of the Laplacian operator is given, up to a constant,  by the difference between $u(x)$ and an average of the neighbors of $u(x)$.  On a non-uniform grid, the operator is given by a similar formula, except the average is replaced by a weighted average (see the discretizations in Section~\ref{sec:3} below.  Each of these discretizations are in the degenerate elliptic form.   In addition, adding a constant term, which corresponds to the inclusion of a term $f(x)$ (which does not depend on $u$ maintains this form.  Furthermore, we can take the maximum or minimum of two terms, and, since the max and min functions are non-decreasing in their arguments, this type of nonlinearity is still degenerate elliptic.

For the most of the discretizations we present below, we can use a clever combination of the Laplacian, and the maximum or minimum terms to produce a discretization which is degenerate elliptic.  This means that we can appeal to the convergence theory for numerical schemes set out in the references above to ensure that the methods converge.  The main new step in the discretization which has not be used in the preceding reference is the use of the irregular grids, which is explained in detail below.

\end{remark}

\subsection{Boundary conditions and far field boundary conditions}
Here we show how to include boundary conditions on more general (non-rectangular) domains, as well as far field boundary conditions.  These boundary conditions are combined along with the elliptic PDE operator into a single (possibly discontinuous) elliptic operator, which is combined with a refinement criteria to perform the grid adaptation.

\begin{example}\label{ex:Poisson}
Consider the Poisson equation
$-\Delta u = f$, with $f$ supported on the unit ball and the far field boundary condition
$u \to 0$, as $\norm{x} \to \infty$.  
An adaptive grid allows us to capture fine details for $x$ near $0$ while reducing computational effort in the far field. 
The artificial boundary condition,
\[
\partial_r u + \frac{u}{r} \approx 0, \quad \text{ for } r \gg 0,
\]
approximates the solution with accuracy ${\mathcal O}\left(\frac{1}{r^3}\right)$~\cite{bayliss1982boundary}. 
\end{example}

\begin{remark}[Characteristic functions on adaptive grids]
Given a domain $\Omega \subset \R^2$, defined the characteristic function of the set $\Omega$ by 
\[
\chi_{\Omega}(x) = 
\begin{cases}
1 & \text{ if } x\in \Omega \\
0 & \text{ otherwise}.
\end{cases}
\]
The characteristic function on a uniform grid was used in \cite{ObermanSINUM} to give a very coarse representation of a boundary (or free boundary).  This representation leads to a piece-wise linear approximation of the boundary of the domain, by connecting  the boundary grid points.   On the adaptive grid, the we obtain a piece-wise linear approximation of the boundary, at difference grid scales, corresponding the spacing of the local grid points.
\end{remark}

\begin{example}\label{ex:Dirichlet}
Consider the Dirichlet problem for the domain $\Omega \subset B = [0,1]^2$ in $\R^2$,
\[
- \Delta u(x) = f(x), \quad \text{ for } x\in \Omega
\]
along with boundary conditions
\[
u = g, \quad \text{ for $x$ on } \partial \Omega
\]
Define the  operator,
\bq\label{Fbc}
F^{bc}[u] = \chi_{\Omega}(x) \left(-\Delta u(x)+f(x)\right) + \chi_{\Omega^c}(x) (u(x)-g(x)) 
\eq
where $\chi_S$ is the characteristic function of the set $S$.  Note that since the characteristic functions are non-negative, this equation is of the form \eqref{GFbc}, so the operator is degenerate elliptic.  Note also that the operator is discontinuous in $x$.

In Section~\ref{sec:3} below, we present the discretization of the Laplace operator on the grid.  In all cases, the approximation has the property that  
\bq\label{discLap}
-\Delta u(x_i) = \bar w u(x_i) - \sum_j w_j u(x_j), 
\quad \text{ where $\bar w = \sum_j w_j$ and $w_j \ge 0$},
\eq
where $u(x_j)$ represents the neighbors of $x_i$.  (On a uniform grid, each $w_j$ would be  equal to  $1/h^2$, where $h$ is the grid spacing.)
This leads to a discretization of \eqref{Fbc} in the form
\[
c(x) \left(\bar w u(x_i) - \sum_j w_j u(x_j) +f(x_i)\right) + d(x) (u(x_i)-g(x_i)) 
\]
with $c(x), d(x) \ge 0$. It is degenerate elliptic according to Definition~\ref{defn:degenElliptic}.  

\end{example}

We can  impose other (for example Neumann or Robin boundary conditions), by replacing the second term with
\[
\chi_{\Omega^c}(x) H(Du(x), x) = 0,
\]
where $H(Du(x), x)$ is itself a first order degenerate elliptic operator.

\subsection{Including free boundaries  in a single degenerate elliptic operator}

The obstacle problem can be formulated as a variational inequality \cite{kinderlehrer2000introduction, glowinski1984numerical}, which is naturally discretized using finite element methods \cite{brezzi1977error}, and solved using a multigrid method~\cite{graser2009multigrid}.

Our approach of adaptive finite difference methods is natural for the obstacle problem, using a formulation of the problem as a degenerate elliptic PDE, however there are far fewer works which use this approach.  Our framework leads to a simple, effective finite difference method which achieves good results using adaptive grids.

We describe here  how to write a free boundary problem as a single degenerate elliptic operator. 

\begin{example}\label{ex:obstacle}
The obstacle problem, for a given obstacle function $g(x)$, which requires that $u(x) \ge g(x)$ and that 
\[
-\Delta u(x) = 0, \quad \text{ for $x$ in  }  \{u(x) > g(x)\}
\]
can be written as a single elliptic equation
\bq\label{obs1}
F^{obs}[u] = \min( - \Delta u, u - g ) = 0
\eq
As in the preceding example, we can discretization the Laplacian in the form \eqref{discLap} and obtain an equation of the form
\bq\label{eq:Obs3}
\min \left (
\bar w u(x_i) - \sum_j w_j u(x_j), u(x_i) - g(x_i)
\right ) = 0,
\eq
which is degenerate elliptic according to Definition~\ref{defn:degenElliptic}.  
\end{example}

This example generalizes to double obstacle problems (using a maximum as well as a minimum), as well as obstacle problems involving nonlinear PDEs which replace the Laplacian~\cite{ObermanSINUM}.

\begin{example}\label{ex:Stef}
The evolution of the one-phase Stefan problem in two dimensions,
$$
\begin{cases}
u_t - \Delta u = 0 &\text{ in } \{u > 0\}\\
u_t - |Du|^2 = 0  &\text{ on } \partial\{u=0\}
\end{cases}
$$
can be represented by the  degenerate elliptic operator,  
\[
u_t + F^{Stef}[u] = 0
\]
with 
\[
F^{Stef}[u] =
\begin{cases}
-\Delta u, 					  &\text{ in } \{ u > 0 \}\\
\min(-\Delta u, - \abs{Du}^2 )          & \text{ in } \{ u \leq 0 \}.
\end{cases}
\] 
More general one phase free boundary problems can be represented as a single operator on the extended domain, as in~\cite{ObermanSINUM}.
Here we have extended from the free boundary to a larger domain, and we solve for the extended operator in the whole domain. 
\end{example}

\subsection{Comparison and Stability of degenerate elliptic finite differences methods}
Stability of degenerate elliptic equations is demonstrated in settings in the reference \cite{ObermanSINUM}.  First, it is shown that there is an explicitly calculated time step so that the forward Euler method is a contraction in the maximum norm.  Second, it is shown that the equation satisfies a nonlinear comparison principle.  While the proof is detailed, some intuition for the results can come from the fact that we can regard the elliptic finite difference equation as expressing a nonlinear average, a point of view taken explicitly in \cite{manfredi2012nonlinear}.  We give a heuristic explanation of these ideas in this remark and refer the reader to the two cited references for more details.

First notice that it is too much to ask that our numerical schemes satisfy the maximum principle.   Even for the equation $-\Delta u = f$ in $\Omega$  with Dirichlet boundary conditions $u=g$ on $\partial \Omega$, the maximum principle does not hold, unless we assume $f \le 0$.   The comparison principle takes the general form 
\[
F^{bc}[u_1](x) \le F^{bc}[u_2](x) \text{ for all $x$} \quad \implies \quad u_1 \le u_2
\]
where the notation $h_1 \le h_2$ means $h_1(x) \le h_2(x)$ for all $x$ (the  domain of definition of the functions is implicit).  
A more specific, and also more explicit form comes from writing $u_1 = S(f_1, g_1)$ and $u_2 = S(f_2, g_2)$ for the solutions of the equation with data $f = f_1, g=g_1$ or $f_2, g_2$, respectively.  In this case, the comparison principle becomes 
\[
f_1  \le f_2 \text{ and } g_1 \le g_2   \quad \implies \quad u_1 \le u_2, 
\]
where $u_1 = S(f_1,g_1), u_2 = S(f_2,g_2)$. 

The discrete comparison principle holds for the numerical scheme, provided that the Laplacian is discretized using an elliptic scheme.  This principle  can be proved using the degenerate elliptic property for a general class of equations which satisfy mild analytical conditions on the operator~\cite{ObermanSINUM} or directly from specific classes of equations without assuming analytical conditions~\cite{manfredi2012nonlinear}.  Once the comparison principle is established, uniqueness of solutions of the schemes follows, since if $f_1=f_2$ and $g_1=g_2$, then $u_1\le u_2$ and also $u_1 \ge u_2$.

The actual proof of the comparison principle, as in the PDE setting, is a proof by contradiction. 
However, we include a plausibility argument which is gives an idea of the reason the local condition can lead to comparison, because it is instructive.  Fix $g_1= g_2$. Starting from \eqref{discLap} and solving for the reference variable, we obtain 
\bq\label{ee1}
u(x_i) =  \sum_j \frac{w_j}{\bar w}  u(x_j) + \frac{f(x_i)}{\bar w}.
\eq
where $x_j$ represents neighbors of $x_i$, and $\sum_j w_j = \bar w$.  This last equation expresses the fact that $u(x_i)$ is a weighted average of its neighbors, plus a constant proportional to $f(x_i)$.  
From this form of the equation, it is plausible that increasing $f(x_i)$ does not decrease the value $u(x_i)$, which leads to the comparison principle.   

The argument which leads to a comparison principle does not depend on the fact that the equation is linear.  So we make a parallel argument for a nonlinear elliptic PDE, with a elliptic finite difference discretization.  
Consider the example of the obstacle problem \eqref{obs1}.  Again the general comparison principle takes the form 
\[
F^{obs}[u_1] \le F^{obs}[u_2]  \quad \implies \quad u_1 \le u_2
\]
which we write in the explicit form
\[
g_1 \le g_2 \quad \implies \quad u_1 \le u_2
\]
where $u_1 = S(g_1), u_2 = S(g_2)$. 

The degenerate elliptic discretization was given in \eqref{eq:Obs3}.  Again, we will solve for the reference variable.  Since we are looking for a solution, which corresponds to the right hand side equal to zero, 
we can divide  each side of the equation by $\bar w = \bar w_i$.  Multiplying the second equation by $\bar w$ also  does not change the solution.  Then we can pull out the term $u(x_i)$ from the equation, and obtain the following equation
\[
u(x_i) = \max \left ( 
\sum_j w_j u(x_j), g(x_i)
\right ) 
\]
Whereas in the linear case, \eqref{ee1}, $u(x_i)$ was an affine function of an average of the values of $u$ at the neighbors, and the data, 
 now $u(x_i)$ is the maximum of that same  average and the data.  In both cases, comparison is suggested by the fact that increasing any of the values of the neighbors or that data does not decrease $u(x_i)$, and stability is suggested by the fact that increasing the values of just one of the neighbors of $x_i$ can increase $u(x_i)$ by at most that amount (and no more).

\section{Adaptive grid}
Our adaptive grid is implemented using a quadtree representation \cite[Chapter 14: Quadtrees]{de2000computational}.  Conceptually, the domain is divided into rectangular regions such that the side length of each neighboring rectangle is either twice, half, or the same as its neighbor. The collection of all vertices are the grid nodes for computing the unknown function. Internally the quadtree is represented as a sparse matrix where the indices of non-zero entries represent coordinates on a fixed ultra-fine grid. 

Our implementation with sparse matrices in MATLAB.
The tool is modular, and the inputs are simple:  the discretization of the operator, $F$, and an additional operator, $G$, used as the refinement criteria, which can be intrinsic (simply setting $G = F$), or defined by the user.   In addition,  if Newton's method is to be used as a solver, the formal Jacobian of the operator is needed, $DF$.  

\subsection{Quadtree construction}
A quadtree is uniquely determined by a list of coordinates and corresponding maximum length scales. Either there is a node at each coordinate with all neighbors within the specified distance, or the coordinate must lie within a rectangle no larger than indicated.  

\begin{figure}
\begin{center}
\scalebox{.25}{\includegraphics{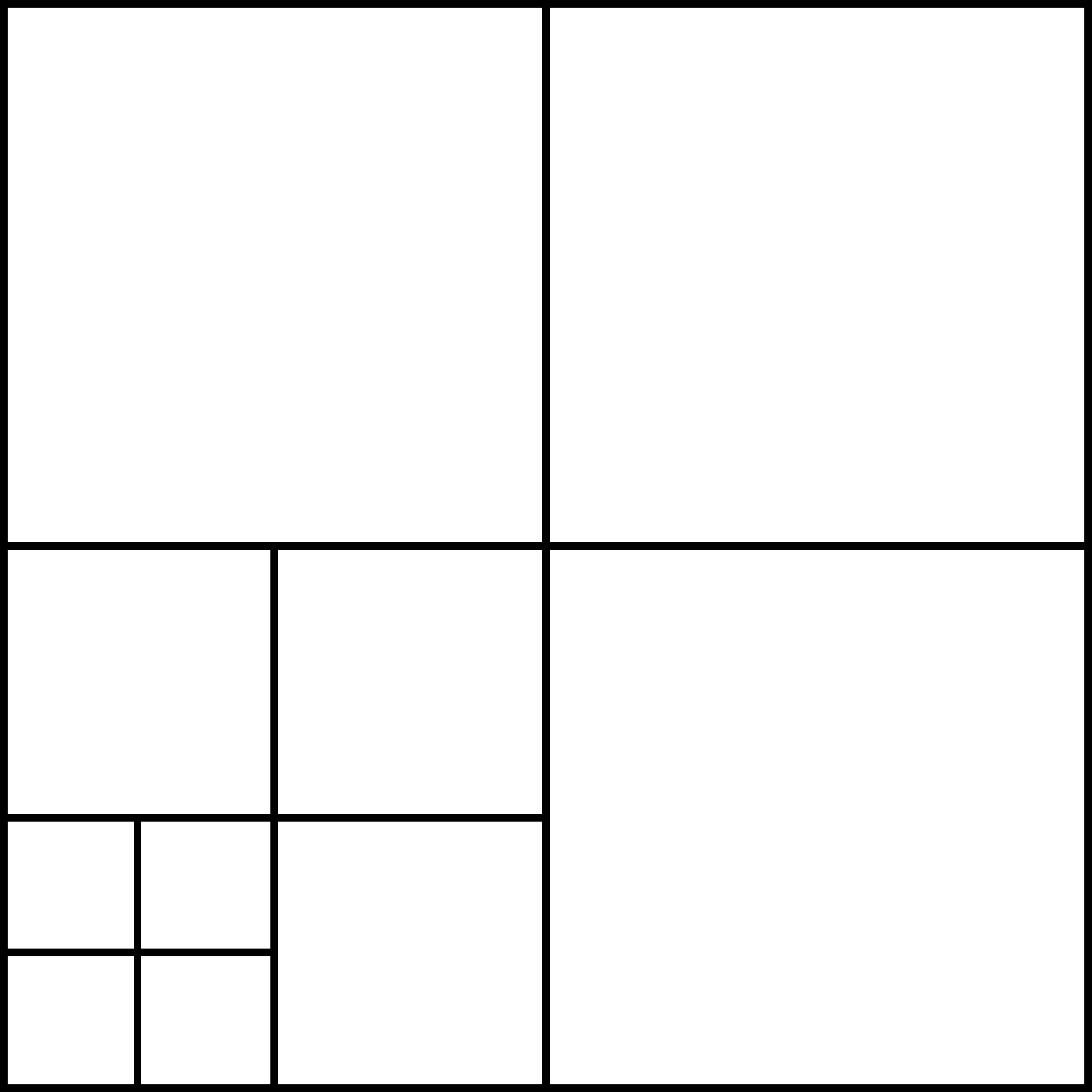}}\hspace{2mm}\scalebox{.25}{\includegraphics{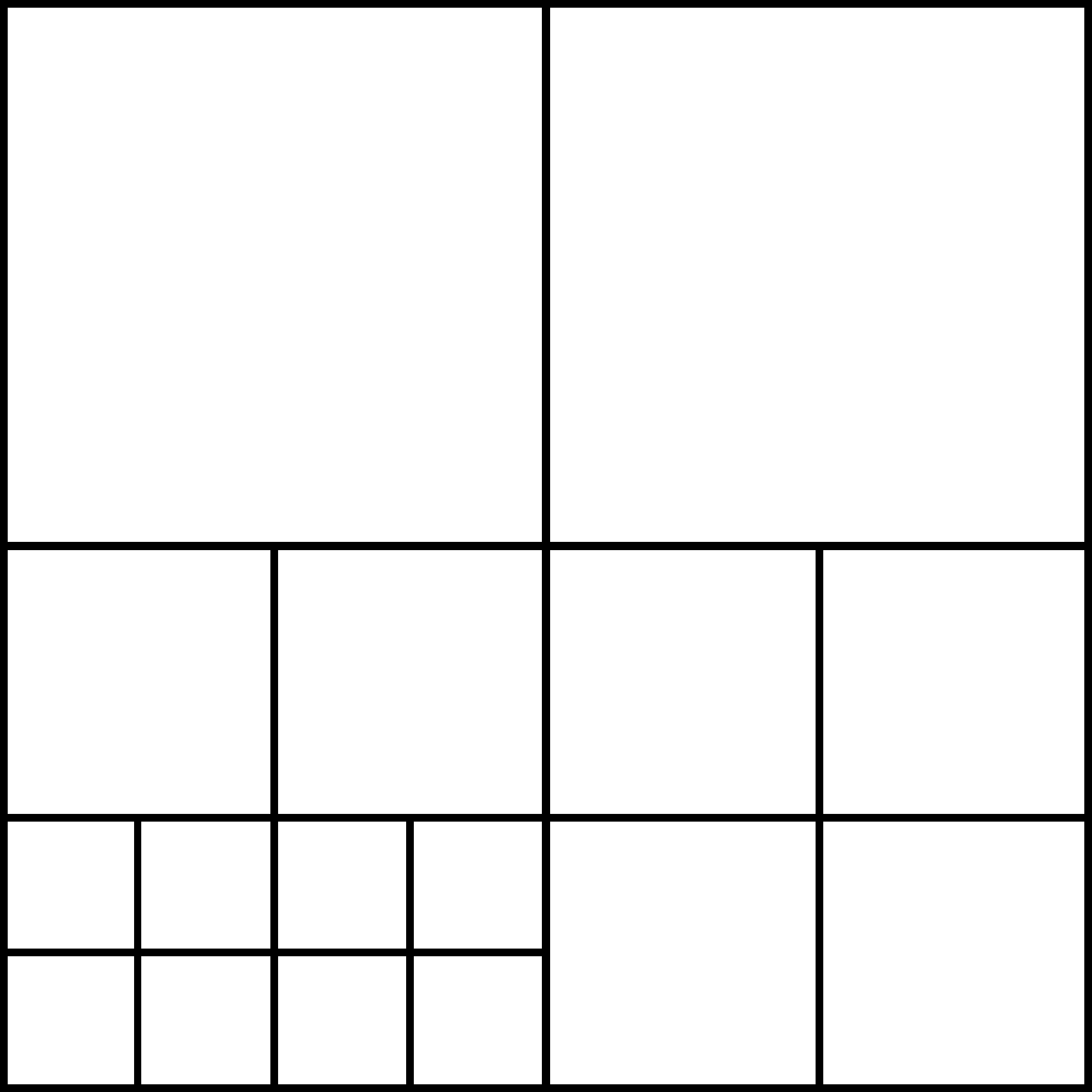}}
\end{center}
\caption{(a) Sample quadtree, consistent with scale-padding constraints $\pad_x=1$ and $\pad_y = 1$.
 (b) Refinement of (a) to be consistent with $\pad_x=2$}
\label{Fig-XScaleBad}
\label{Fig-XScaleGood}
\end{figure}

To discretize the Laplacian, we impose an additional `scale-padding' constraints depending on the aspect ratio of the physical domain.  Dangling nodes, vertices with three neighbors, occur midway along the shared edge of two equal-sized rectangles one of which is subdivided.  The scale-padding constraint in $x$ specifies the minimum number of equal-sized rectangles that must exist to both the left and right of a dangling node.  Figure~\ref{Fig-XScaleBad} illustrates a pair of quadtrees, the latter refined to observe the scale-padding constraints $2$ in $x$ and $1$ in $y$.  

To built a quadtree over a virtual ultra-fine grid of $2^N+1$ by $2^N+1$, we build a list of squares the quadtree must contain.
\begin{enumerate}
	\item List all requested or required squares size $2^k+1$ by $2^k+1$.
	\item Add siblings to the list. Fill to the edge, if needed to prevent a dangling node too close to the boundary.
	\item List all parents.  Expand list to satisfy scale-padding constraints.
	\item Use grandparents to ensure the expanded list of parents includes all of their siblings.  
	          This is the list of all required squares size $2^{k+1}+1$ by $2^{k+1}+1$.  
\end{enumerate}
Given $M$ requested coordinate-length scale pairs, this procedure is ${\mathcal O}\left(N M \log M\right)$.

\subsection{Adaptivity}
Any refinement scheme based on position and the local value of the function and its derivatives may be specified. 
A ``refinement criteria'' operator is computed at all current grid points and compared with a ``refinement threshold''.  
The new grid must be refined to the finest available spacing at all nodes where the criteria exceeds the threshold. 

For better control, the user may supply a non-decreasing array of refinement threshold values.
Where the criteria exceeds the $k^{th}$ threshold the new grid is refined to at least the $k^{th}$-finest scale.

The user may also supply a padding parameter, above that required for discretization of the Laplacian at dangling nodes.  It is often convenient to give a simple refinement criteria and a large padding parameter.

All adaptive grids must include the fixed initial quadtree, which is determined by the placement of the initial data.
The initial data is treated as scattered and linearly interpolated onto the smallest quadtree for which every supplied data-point appears as a node.
We assume the placement of initial data implies the minimally acceptable spatial resolution.

\section{Discretization on the adaptive grid}\label{sec:3}

As mentioned above, degenerate elliptic schemes are easily built from the upwind schemes for derivatives and Laplacian.  Our solvers expect the user to specify their operators in terms of these building blocks. For convenience, the discretization discussed here is kept in a black box.

We will consider regular nodes, dangling nodes, and boundary nodes separately.  Dangling nodes are those with only three neighbors and occur midway along the edge of a rectangle that adjoins two half-size rectangles. 

\begin{figure}
\includegraphics[width=0.4\textwidth]{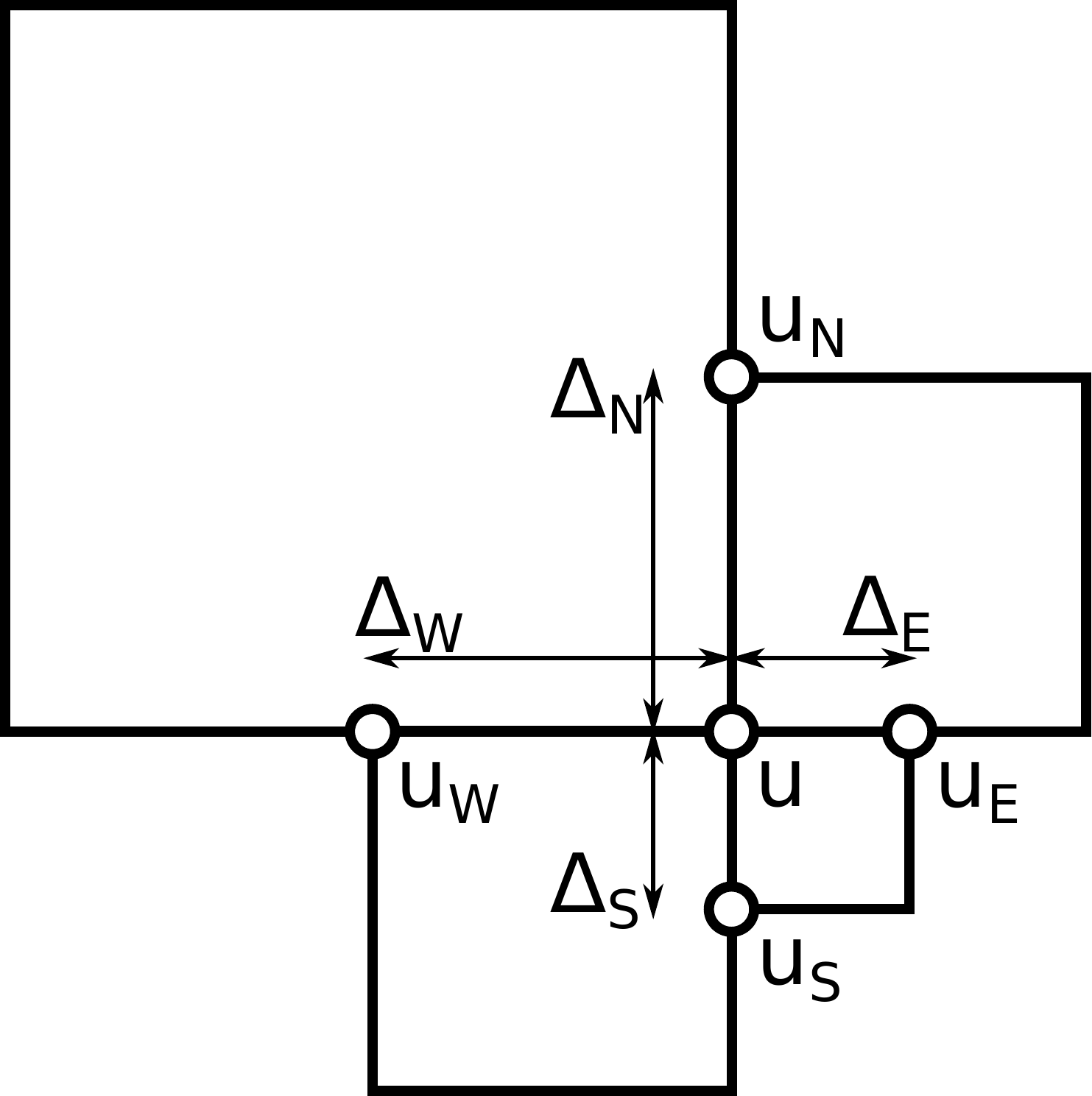}
\hspace{2mm} {\includegraphics[width=0.4\textwidth]{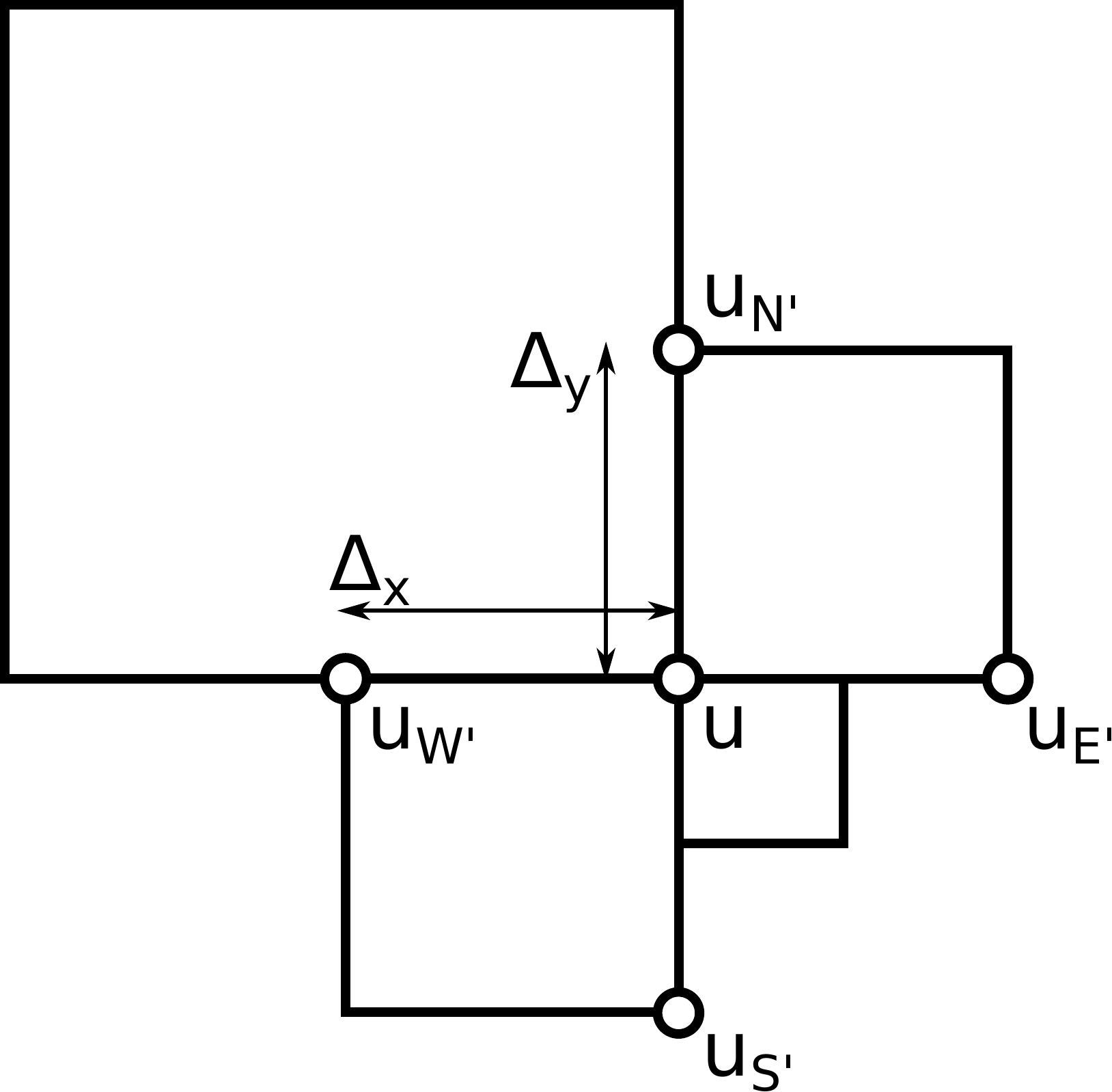}}
\caption{One of six possible configurations at a regular node.  (a)  The stencil for the discretization of first-derivatives. (b) The stencil for the Laplacian discretization.}
\label{Fig-RegularNode} 
\end{figure}

A regular node is the shared vertex of four rectangles. Consider the nearest neighbors $u_E, u_W, u_N$ and $u_S$ at distances $\Delta_E, \Delta_W, \Delta_N$ and $\Delta_S$ respectively, as in Figure~\ref{Fig-RegularNode}.
The standard upwind discretizations are:
\begin{equation}\label{Eqn-DiscDerivReg}\begin{aligned}
\partial_{x}u \approx \frac{u-u_W}{\Delta_W}
&&\text{ and }&&
-\partial_{x}u \approx \frac{u-u_E}{\Delta_E},
\end{aligned}\end{equation}
both accurate to first order. 
For the Laplacian operator we identify the nearest pairs of equidistant opposing nodes, $u_{E'}$ and $u_{W'}$, and  $u_{N'}$ and $u_{S'}$, as in Figure~\ref{Fig-RegularNode}.  The standard discretization for $\partial_x^2u$,
\begin{equation}\label{Eqn-DiscLaplacianNaive}
-\partial_x^2u\approx \frac{2u - u_{E'} - u_{W'}}{2(\Delta_x)^2},
\end{equation}
is accurate to second order.  Discretization using only nearest-neighbors is only accurate to first order.

\begin{figure}
\begin{center}
\includegraphics[width=0.4\textwidth]{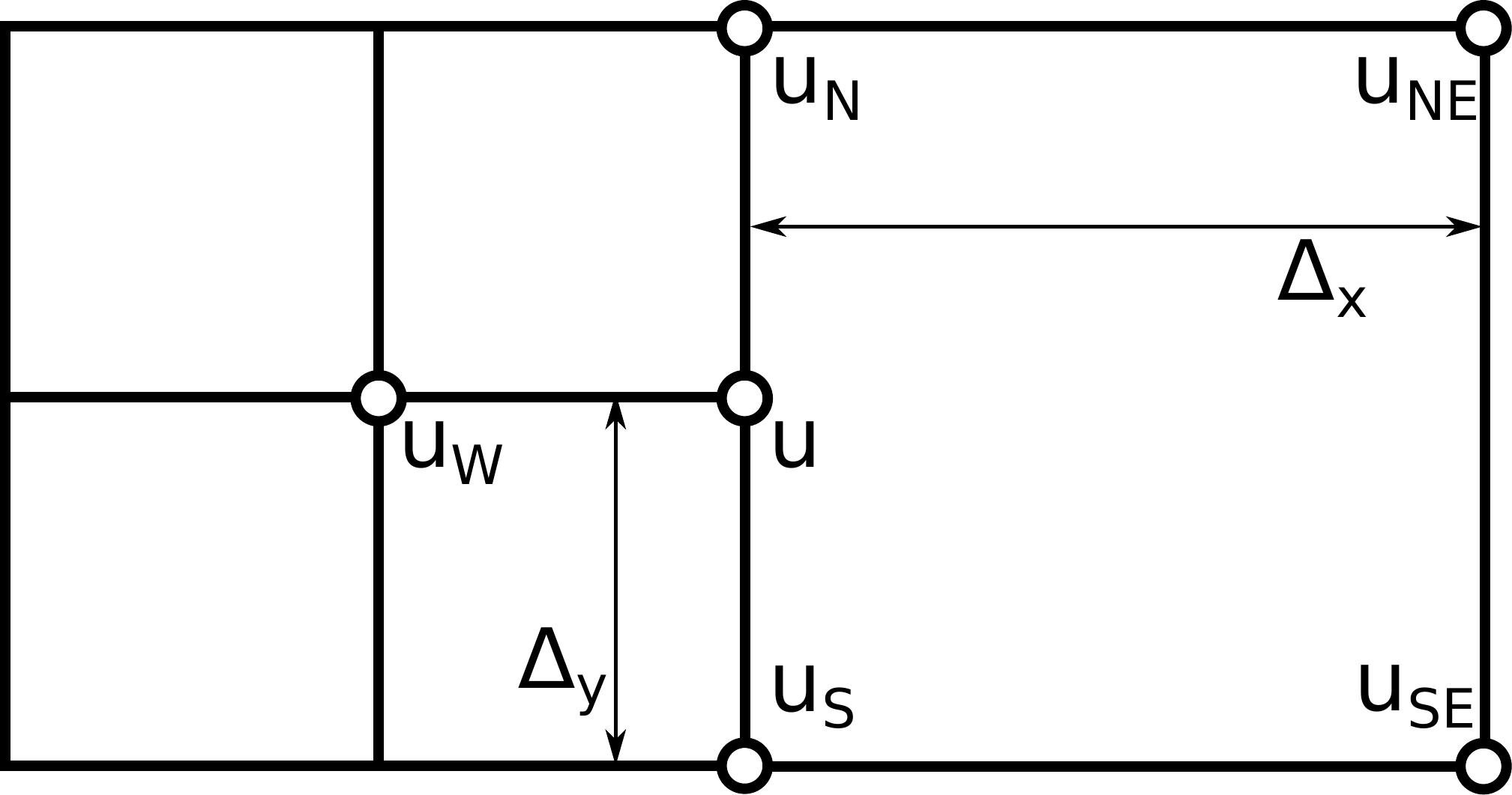} \hspace{2mm} \includegraphics[width=0.4\textwidth]{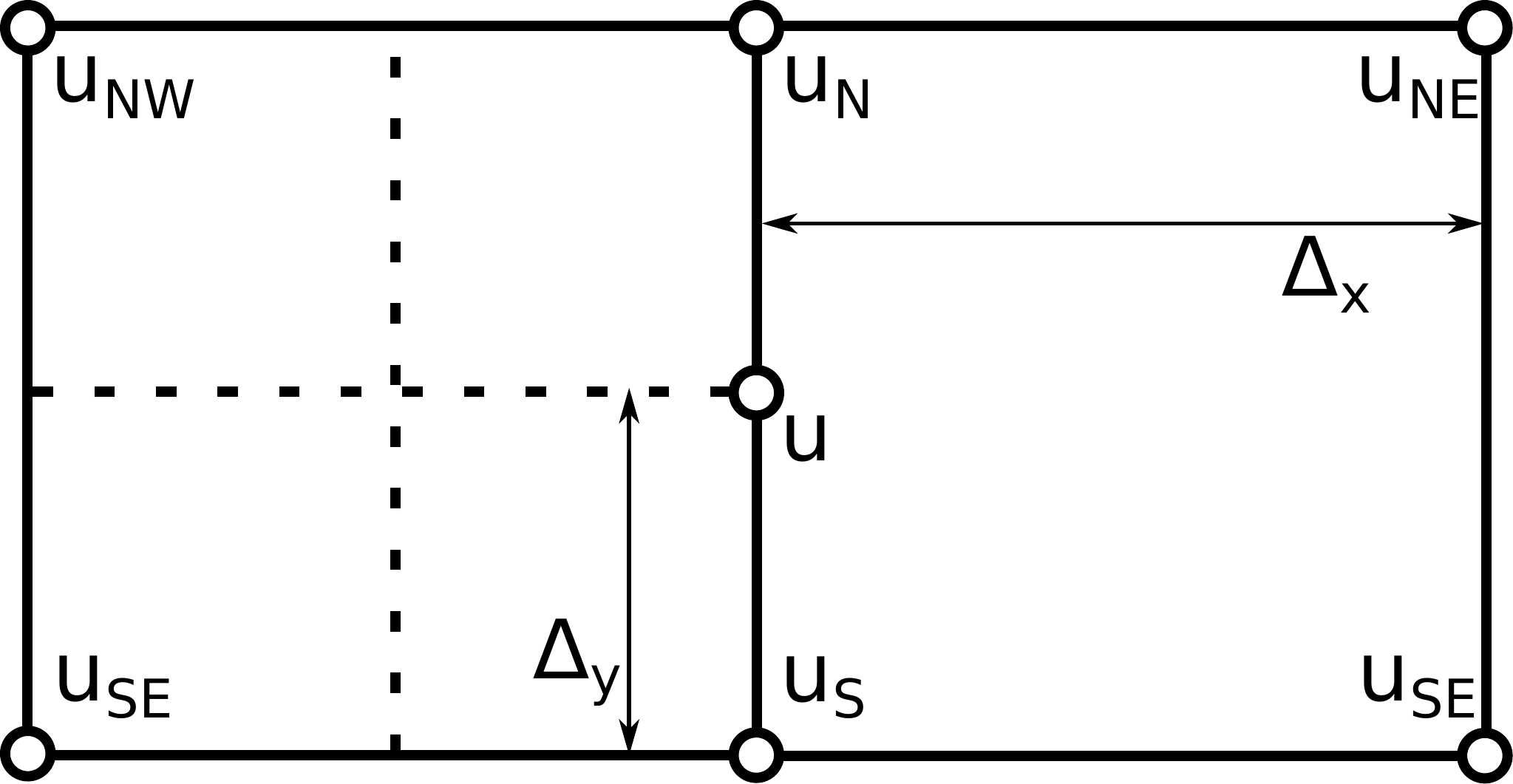}
\end{center}
\caption{(a) A dangling node in the $x$ variable, showing the stencil for first-derivative discretization. (b) The stencil for the Laplacian discretization at a dangling node, provided $\Delta y \leq\Delta x$.}
\label{Fig-DanglingNode}\label{Fig-IStencil-Dangling}
\end{figure}

For a dangling node, as in Figure~\ref{Fig-DanglingNode}, we use the farther vertices of the larger square to interpolate a value for the unknown function directly opposite, then discretize as at a regular node. In the illustrated situation we would use,
\begin{equation}\label{Eqn-DiscDerivDang}
\partial_{x}u \approx \frac{u - \frac{u_{NE}+u_{SE}}{2}}{\Delta_x}.
\end{equation}
This is a monotone discretization, accurate to first order since $u_N$ and $u_S$ are equidistant.
Second derivatives are more difficult.  There is no general upwind discretization for $\partial_{x}^2u$ for the dangling node illustrated in Figure~\ref{Fig-DanglingNode}.

\begin{figure}
\begin{center}
\scalebox{.25}{\includegraphics{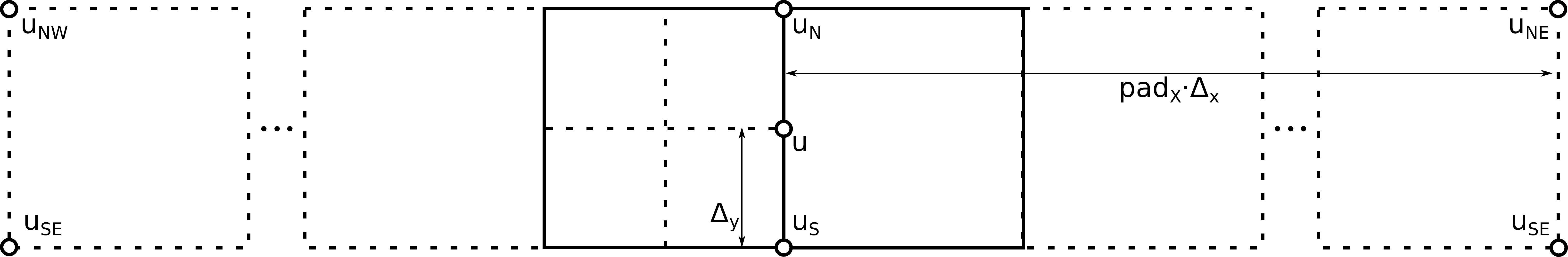}}
\end{center}
\caption{The extended stencil for laplacian discretization at a dangling node, when $\Delta y \not\leq\Delta x$.}
\label{Fig-IStencil-Wide}
\end{figure}

At dangling nodes, we choose to discretize the laplacian using an I-shaped stencil as in Figure~\ref{Fig-IStencil-Dangling}.  The following expansions are accurate to second-order:
\[\begin{aligned}
2u - \frac{u_{NW}+u_{SW}+u_{NE}+u_{SE}}{2} &= -\left(\Delta y\right)^2\partial_y^2 u -\left(\Delta x\right)^2\partial_x^2 u + {\mathcal O}(\Delta x^4 + \Delta y^4)\\
2u - \left(u_N+u_S\right) &= -\left(\Delta y\right)^2\partial_y^2u + {\mathcal O}(\Delta y^4)
\end{aligned}\]
Our discretized Laplacian is then:
\begin{equation}\label{Eqn-DiscLaplacian}\begin{aligned}
-\partial_y^2u - \partial_x^2u \approx &\frac{1}{(\Delta x)^2}\left(2u - \frac{u_{NW}+u_{SW}+u_{NE}+u_{SE}}{2}\right)\\ &+ \left(\frac{1}{(\Delta y)^2} - \frac{1}{(\Delta x)^2}\right)\left(2u - \left(u_N+u_S\right)\right).
\end{aligned}\end{equation}
This is a monotone discretization provided $\Delta y \leq \Delta x$.  Should $\Delta y > \Delta x$ we take a wider I-shaped stencil, as in Figure~\ref{Fig-IStencil-Wide}. 
 To ensure $\Delta y \leq \pad_x\Delta x$, (or, $\Delta x \leq \pad_y\Delta y$, for the other type of dangling node,) we choose,
\[\begin{aligned}
\pad_x = \left\lceil \frac{L_y}{2L_x} \right\rceil.
&& \text{ and }  && \pad_y = \left\lceil \frac{L_x}{2L_y} \right\rceil.
\end{aligned}\] 
These values reflect the aspect ratio of the domain in physical variables.  The `scale-padding' constraints when building the quadtree guarantee this wider stencil is available at all dangling nodes.

\subsection{Boundary nodes}
At boundary nodes where Dirichlet boundary conditions are not provided, 
we implement generic Robin boundary conditions:
\[
A(x) u'(x) + B(x) u(x) = C(x)
\]
The functions $A$, $B$, $C$ should be provided for each edge of the domain.  
User-friendly shortcuts for Neumann or Dirichlet conditions are provided.  

Where $A(x) = 0$,  $u(x)$ is specified. The node is considered inactive and is not updated by means of a logical mask.
Where $A(x) \neq 0$, we use the boundary conditions to determine the outward derivative and the regular discretization for the inward derivative. We discretize the second derivative by weighting the outward and inward derivatives equally, as in \eqref{Eqn-DiscLaplacianNaive}.


\section{Numerical Solvers}
In this section we discuss the approach to solve the nonlinear  finite dimensional equations obtained by discretization of the elliptic or parabolic PDE on the adaptive grid.

In the case of a nonlinear parabolic (time-dependent) problem, $u_t = F[u]$, it is natural to discretize using the Forward Euler method, which leads to 
\[
\frac{u^{n+1}-u^n}{dt} = F[u^n].
\]
This discretization leads to the explicit iteration
\[
u^{n+1} = u^n + dt F[u^n].
\]
For stability, there is a restriction on the time step $dt$, which is often proportional to $dx^2$ for parabolic problems, where $dx$ is the grid discretization parameter (for uniform grids).   The theory developed in \cite[Theorem 6]{ObermanSINUM} allows us to obtain an explicit value for the time step restriction which ensures stability in the maximum norm.   In addition, we can use local values of this time step, which allows different time steps at different grid resolutions.  In this manner, a global time step can be taken which corresponds to multiple iterations at small scales, and a single iteration at the largest scale.

When the Laplacian is discretized according to \eqref{discLap}, where the weights $w_i$ depend on the local grid, and are equal to $(\Delta_x)^2$ or $(\Delta_y)^2$ in the case of \eqref{Eqn-DiscLaplacianNaive} and in the case of dangling nodes can be obtained from \eqref{Eqn-DiscLaplacian}.   In this case, the restriction on the time step at $x_i$ is given by
\[
dt \le \bar w_i
\]
The local dependence of the time step is encoded directly in the weight $\bar w_i$.  For nonlinear discretizations involving the Laplacian, such as the obstacle problem, given by \eqref{obs1}, the local time step restriction is \emph{the same}.   In general, the local time step restriction is given by the local Lipschitz constant of the scheme, regarded as a function.  So simply differentiating the scheme and taking an upper bound can given an acceptable time step.  In some cases, this constant may depend on the initial data, as in the case of $u_t = \abs{u_x^2}$ where $dt \le dx^2/\max_j {u^0_{j+1}-u^0_j}$ from \cite[Section 4]{ObermanSINUM}.
 
\subsection{Details on the time-dependent solver}
Before time evolution, we apply the refinement scheme to the linear interpolation of the initial data and iterate to ensure infill of any coarse regions of the initial quadtree that are nevertheless of interest.

Nodes are separated according to the distance to their nearest neighbor. A nonlinear CFL condition is required to determine a characteristic time-scale for each group. The time-scales differ by powers of two. We list each group according to the inverse ratio of their characteristic time-scales and the coarsest time-scale. The result is randomly permuted before each time-step to produce a visitation schedule. To evolve by one time-step, all nodes in each group are simultaneously updated according to the visitation schedule.

This scheme is optimal in the sense that each group-update operation is of the order of the number of nodes in the group, and that no more updates occur than required by the nonlinear CFL condition.

\subsection{Static solver}
In the case of elliptic equations, it is possible to find the solution by iterating the time dependent problem.  This follows from the fact that the forward Euler method, with the restricted time step is a contraction,  \cite[Theorem 7]{ObermanSINUM}. 
However, this method is slow, since the time step can scale quadratically with the smallest spatial scale, $dt = \bO(dx^2)$, so the number of iterations to solve to a fixed time grows as the grid grows.

An effective alternative is to implement Newton's method.  We build the exact Jacobian of the discrete scheme, meaning the gradient of the scheme regarded as a function, which requires writing additional code to represent the Jacobian.
In addition, since we are working in two dimensions, and the Jacobian is a sparse matrix, we can use direct solvers effectively.
The Jacobian is sparse with the number of nonzero elements on the order of the number of nodes.  We expect Gaussian elimination in ${\mathcal O}\left(M\log M\right)$ time, where $M$ is the number of nodes.  This results in a fast direct solver.

For further efficiency, we seek a solution at coarse scales before allowing refinement to finer scales. This corresponds to the first step in a V-cycle for a multigrid method.   Starting with the initial quadtree we iterate Newton's method until a stopping condition is reached. We then allow refinement up to the second-coarsest scale present in the initial quadtree as called for by the refinement criteria and threshold. We then seek another solution and repeat the process until all scales are allowed.  

Facility for stopping criteria and thresholds are similar to those for refinement.  When allowing refinement to the $k^{th}$-finest scale, iteration of Newton's method continues until the stopping criteria is less than the $k^{th}$ stopping threshold at all nodes.  The default stopping criteria is the $L^\infty$ norm of the elliptic operator.

Note that when seeking a solution over a fixed multi-scale grid it is more efficient to define the multiscale grid through the refinement criteria and provide initial data only on the  coarse grid.

\section{Computational Examples}
In this section we present numerical results, which show the validity and performance of the method, and allow for us to demonstrate the effectiveness of different refinement strategies.

\subsection{Artificial boundary conditions}\label{ex:artificial}
We are in the setting of Example~\ref{ex:Poisson}:
the Poisson equation $-\Delta u = f$, with $f$ supported on the unit ball and
$u \to 0$, as $\norm{x} \to \infty$.  
Set
\[
f(r,\theta) = r(1-r)^+\sin(5\pi r)\cos(3\theta),
\]
where $(1-r)^+ = \max(1-r,0)$,
and impose the artificial boundary condition,
$\partial_r u + \frac{u}{r} = 0$ at the boundary of the computational domain.
Set the domain to be a square domain with side length $2\times 10^3$.  
The refinement criteria: the grid should be finest for $r<1$.  

\begin{figure}
\begin{center}
\includegraphics[keepaspectratio=false, width=3.5in, height=4in]{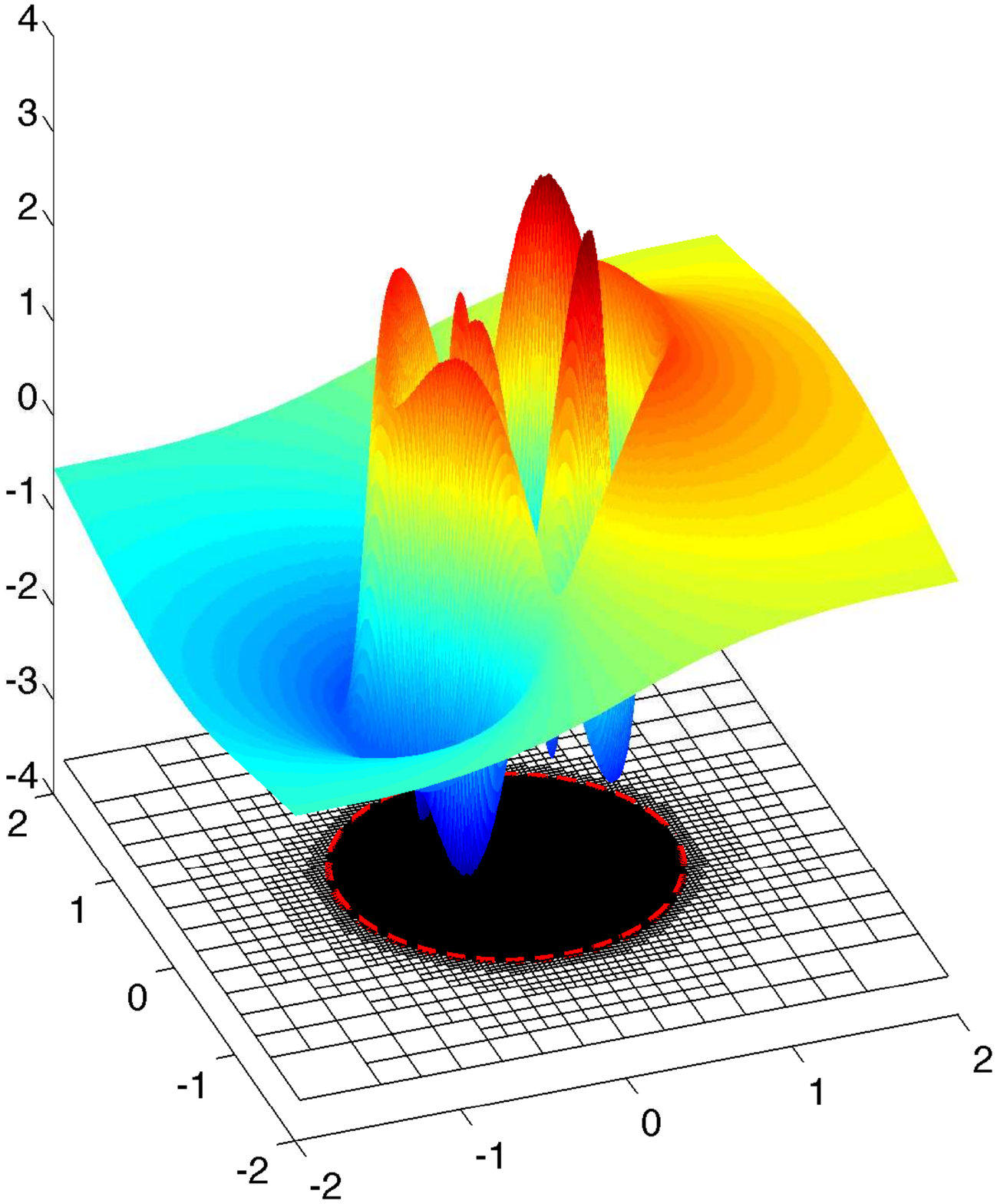}\includegraphics[keepaspectratio=false, width=3.5in, height=4in]{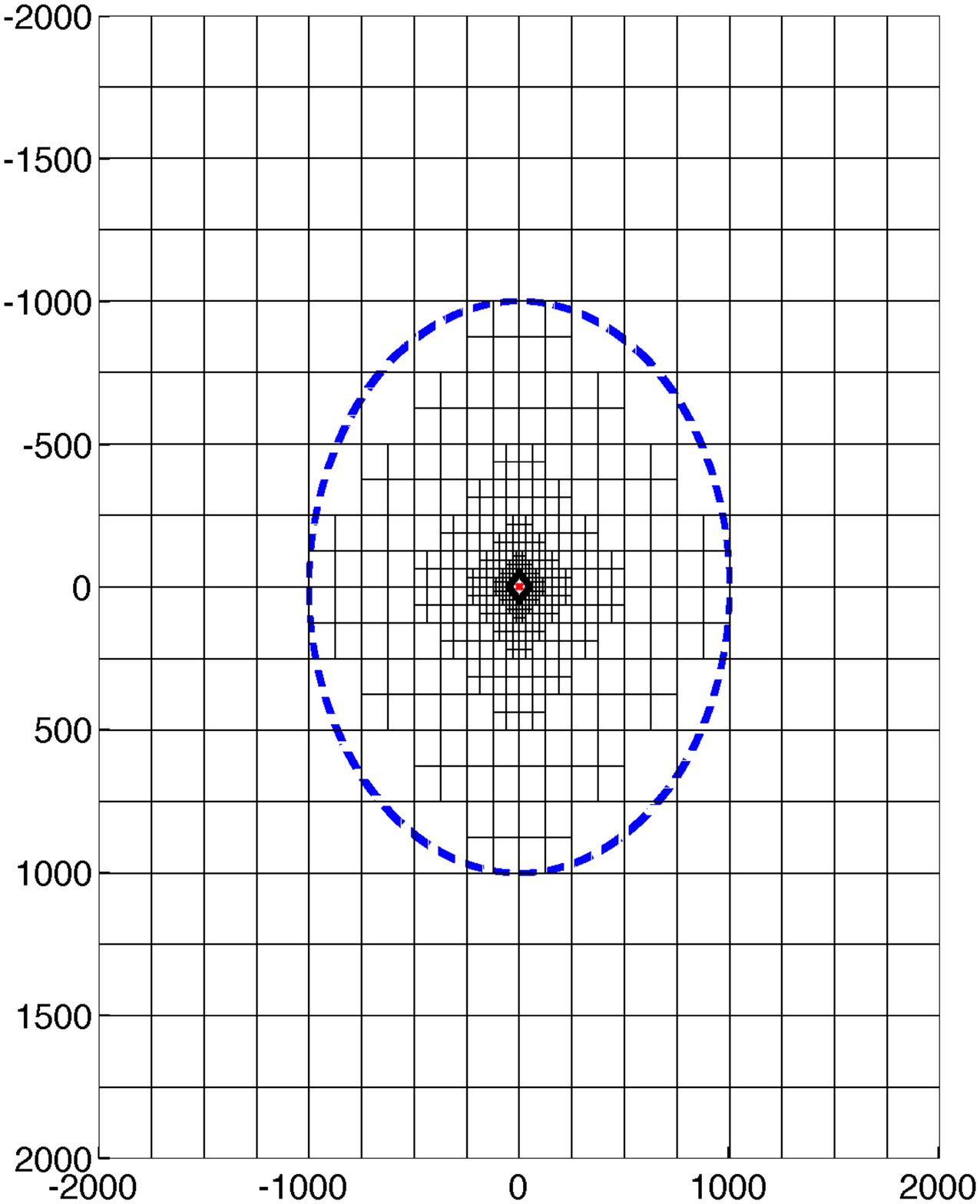}
\caption{
(a) Detail of the solution of the Poisson equation with $f(x)$ supported on for $r\le1$ (indicated by the red dashed circle). (b) The corresponding adapted grid.  Artificial boundary conditions are applied for $r>10^3$ (indicated by the blue dashed circle).}
\label{ex:artificial:figA}
\end{center}
\end{figure}

A detailed view of the solution in the near field can be see in Figure~\ref{ex:artificial:figA}(a). 
The layout of the grid at large scale can be seen in Figure~\ref{ex:artificial:figA}(b).   Table \ref{table:artificial} outlines the allocation of computing resources and nodes by region.     Broadly speaking, the adaptive grid allows us to compute on a very large domain, with a computational cost on the order of (within a couple multiples of)  restricting to the unit square.

\begin{table}
\begin{center}
    \begin{tabular}{c| c| c| c}
	Region & Relative Area & Time Spent & Nodes on Final Grid\\
	 \hline
	 $r < 1$			& $\ll$0.001\%	& 38.5\%	& 78.3\% \\
	 \hline
	 $1<r<10$		& 0.002\%		& 13.3\%	& 15.0\% \\
	 \hline
	 $10<r<10^3$	& 19.6\%		& 26.6\%	& 5.2\% \\
	 \hline
	 $10^3<r$		& 80.4\%		& 21.6\%	& 1.5\%
    \end{tabular}
    \caption{Resource use compared to detail achieved, by domain region, for the example solution of the Poisson equation with artificial boundary conditions.  Artificial Boundary conditions were applied for $r>10^3$.}  
    \label{table:artificial}	
\end{center}	
\end{table}

\subsection{Irregular domains}\label{ex:irregular}
Consider a problem of the type Example~\ref{ex:Dirichlet}, where the Dirichlet problem is posed on an irregular 
domain $\Omega$, contained in a rectangle. 
For the Poisson equation with Dirichlet Boundary conditions we use the operator $F^{bc}$.

\begin{figure}
\begin{center}
\includegraphics[keepaspectratio=false, width=5in, height=5in]{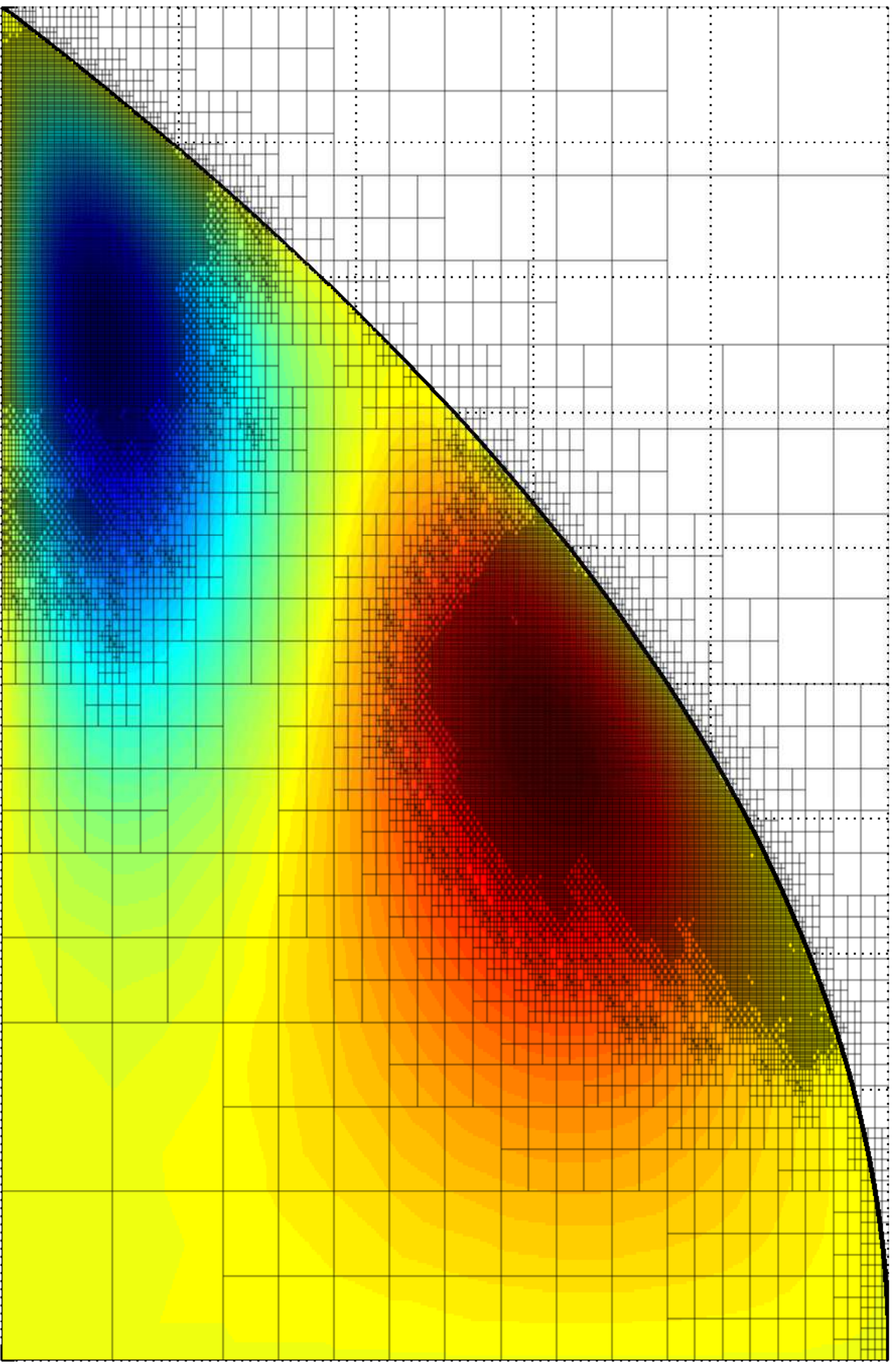}
\caption{Solution of a Poisson equation on a curved domain with Dirichlet boundary conditions.}
\label{ex:irregular:figD3}
\end{center}
\end{figure}

The refinement criterion used was based on the combination of residual of the operator and the proximity to the boundary.    An example of what can be accomplished is shown in Figure~\ref{ex:irregular:figD3}.
Notice that this leads to a maximal refinement in two blobs near the boundary, near local extreme points of the solution (red, where the solution value is near 60 and blue, where the solution value is near -100), while other areas near the boundary have a relatively coarse grid (yellow, where the solution is near 0).

\begin{figure}
\begin{center}
\includegraphics[keepaspectratio=false, width=4.5in, height=4.5in]{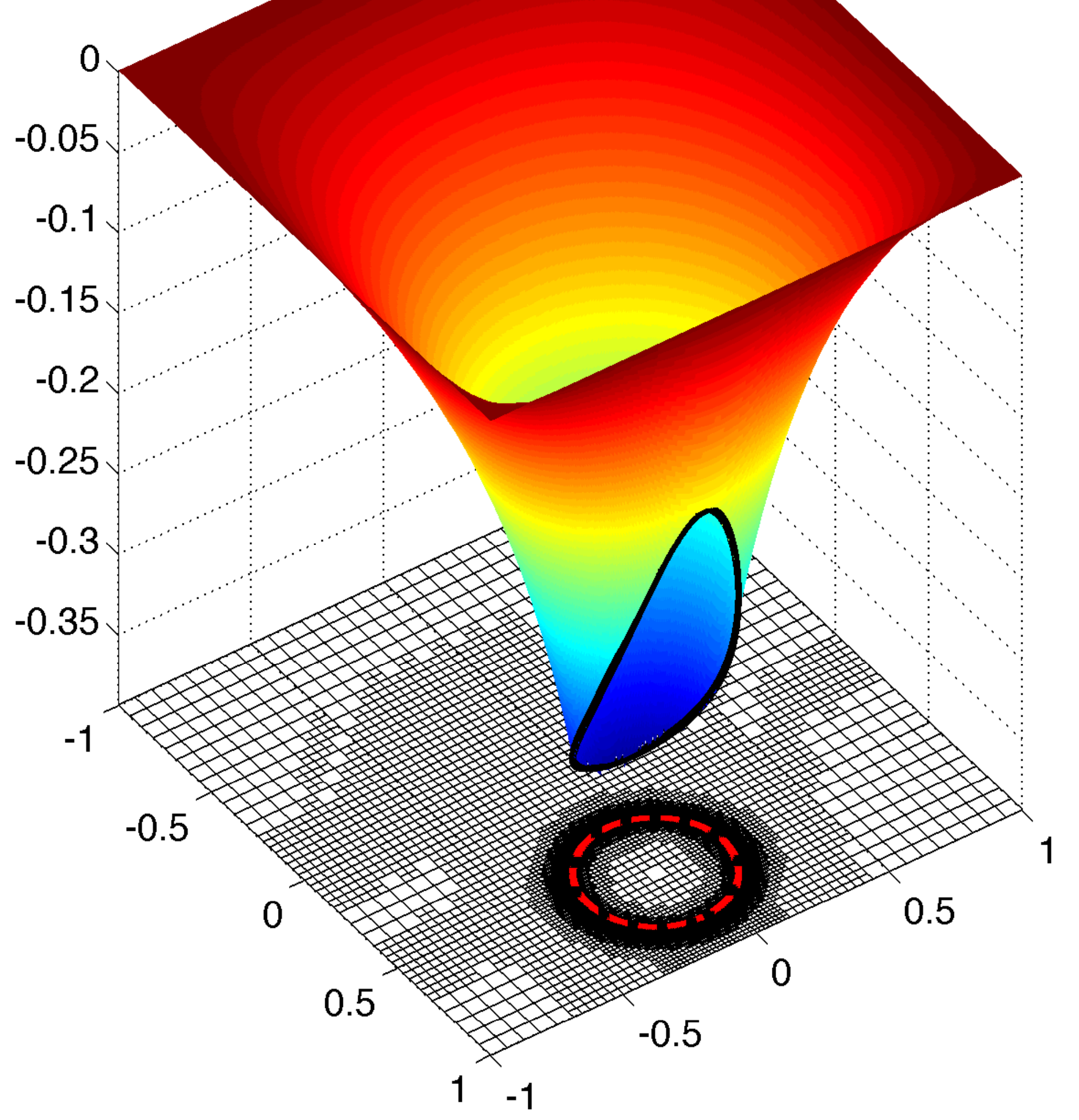}
\caption{Solution of Laplace equation on a punctured domain, with inhomogeneous Neumann boundary conditions on the red circle.}
\label{ex:irregular:lap_figA}
\end{center}
\end{figure}

To impose Neumann or Robin Boundary conditions, we apply the boundary conditions for grid points near the boundary, and further away, simply impose $u=0$.  As an example, Figure \ref{ex:irregular:lap_figA} presents results for the Laplace equation with homogeneous Dirichlet boundary conditions on the boundary of the unit square, combined with inhomogeneous Neumann  boundary conditions 
\[
\frac{ du }{ dn} = 1, \quad \text{ for $x$ on the boundary of a punctured circle inside the domain}.
\]
The adaptive grid was determined by the slope of the solution combined with proximity to the interior boundary.

\subsection{Obstacle problems}\label{ex:Obstacle}
We are in the setting of the obstacle problem, Example~\ref{ex:obstacle}, represented by the operator $F^{obs}$.

\begin{figure}
\begin{center}
\includegraphics[width = 5in, height = 4in]{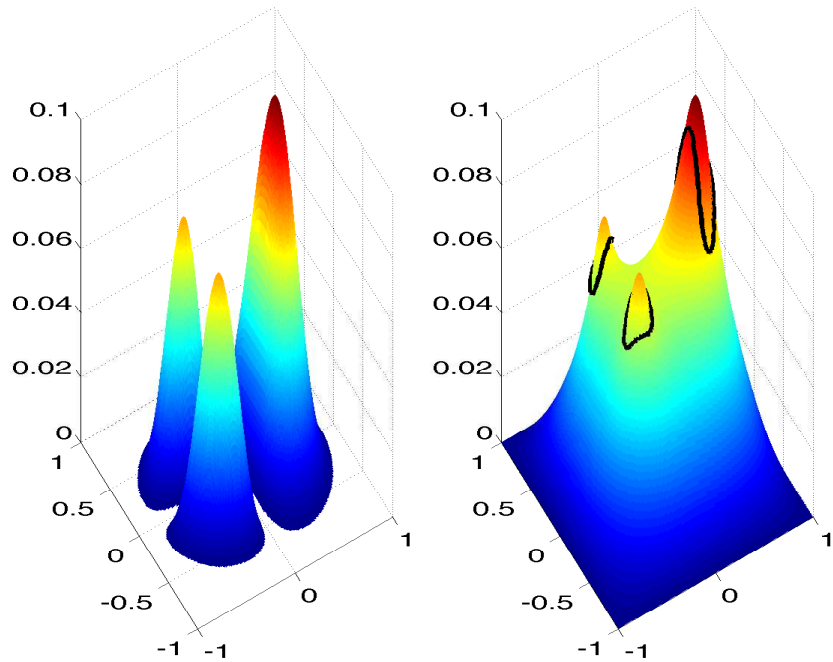}
\caption{(left) Obstacle. (right) solution of the obstacle problem (with contact contour in black).}
\label{ex:obstacle:figA6}
\end{center}
\end{figure}

The obstacle $g(r,\theta)=r^2\cos^2(\theta)$, multiplied by a factor of $2\sin^2(\pi y)$ for $x<0$, and by $\exp(-r)$ for $r>\frac{1}{4}$. The obstacle and solution are shown in Figure~\ref{ex:obstacle:figA6}. 
 
The problem was solved using different refinement criteria, defined as follows. The contact contour determined in all three cases was virtually indistinguishable.

 - As a baseline method, we used a simple predetermined (non-adaptive) grid criteria.  The finest grid resolution is specified by the distance to the local maxima $x_1, x_2, x_3$ of the obstacle (since the contact set is unknown).
 \[
 G(x)  \text{ determined by }( \text{dist}(x, \{ x_1, x_2, x_3 \} ))
 \]
In this case, the solution was found using nonlinear multigrid, allowing a progressively finer grid each time the residual drops below a threshold.  

\begin{figure}
\begin{center}
\includegraphics[height = 5 in]{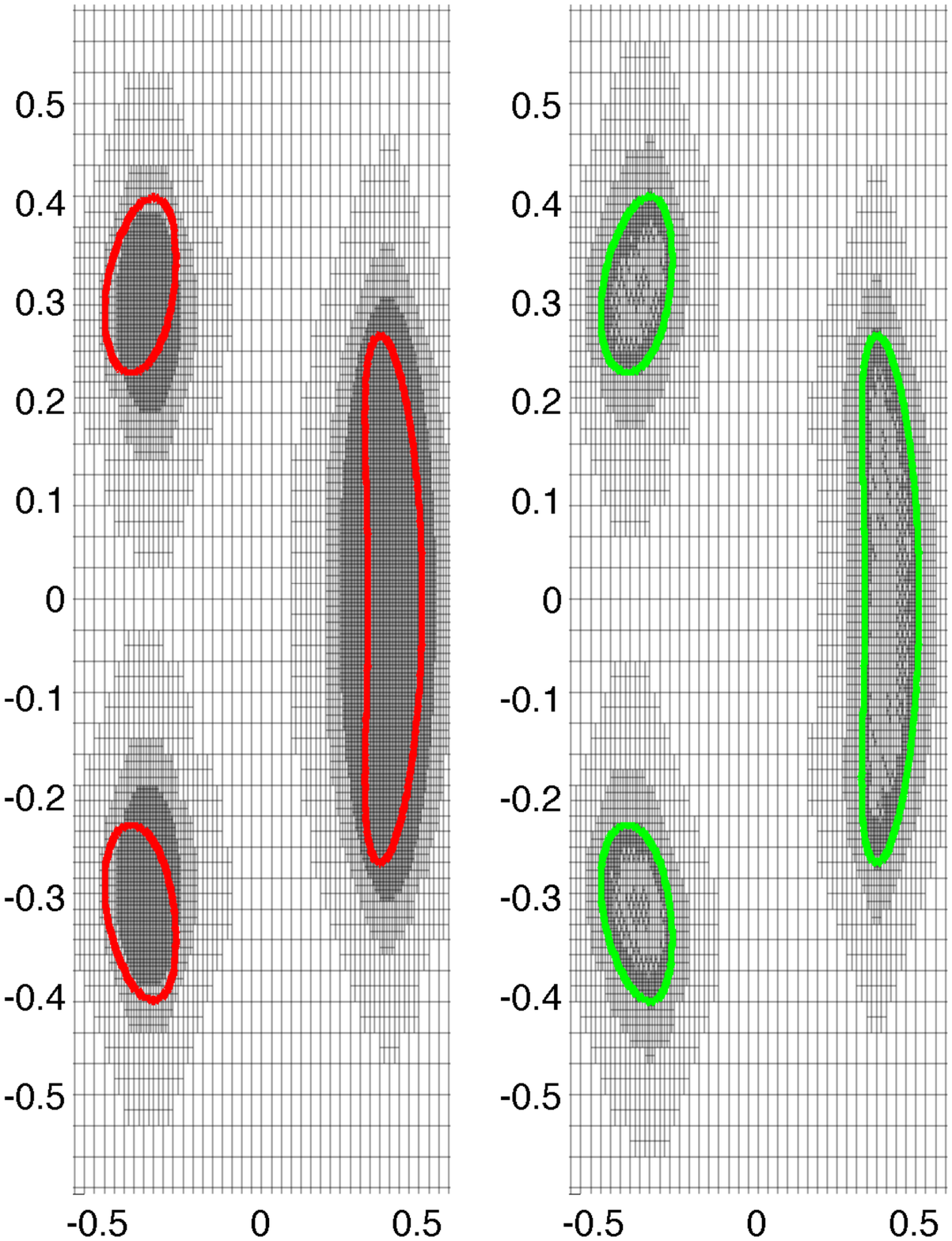}\hspace{2mm}\includegraphics[width= 4. in]{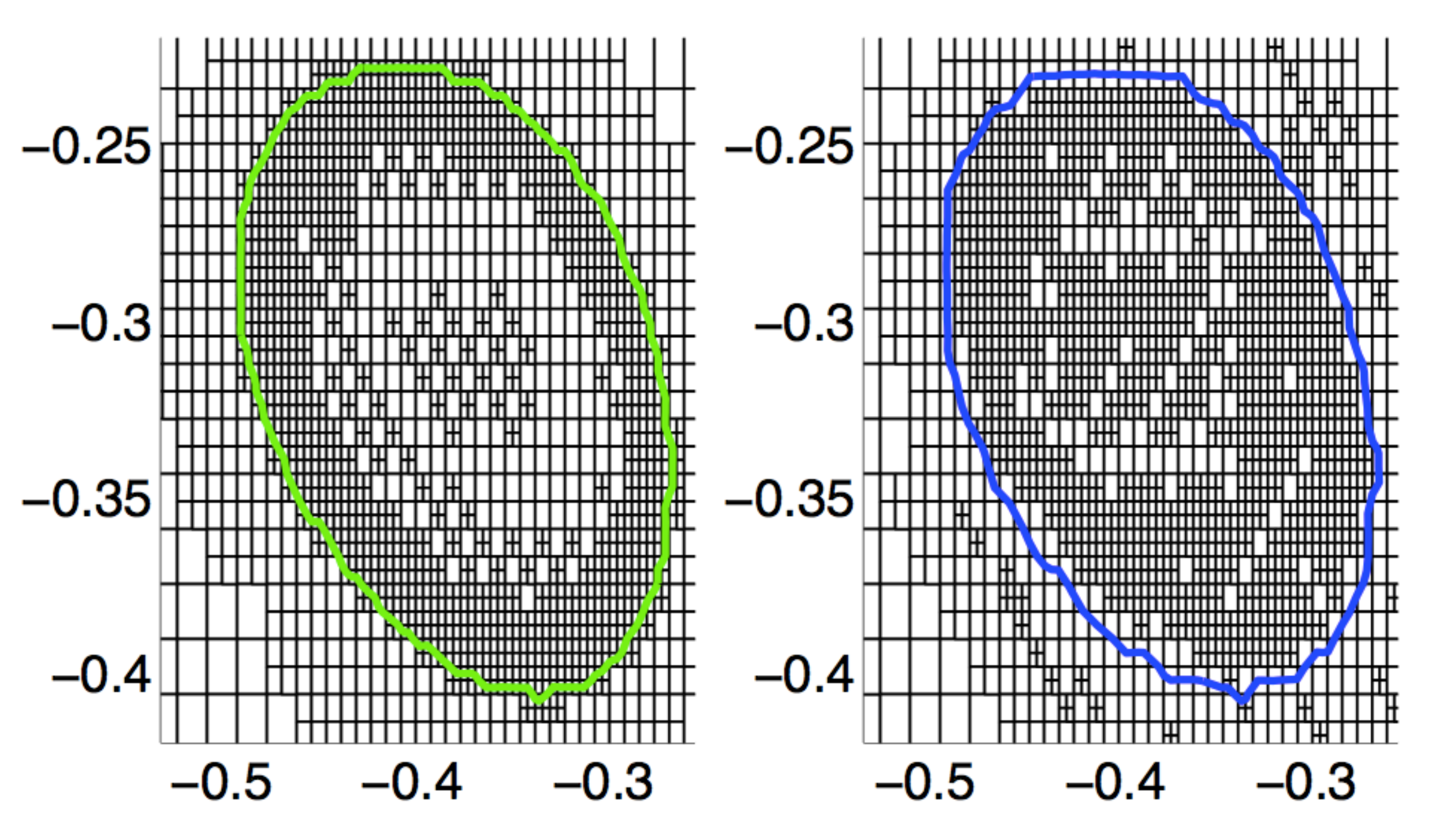}
\caption{(top) Detail of predetermined grid (left) and boundary-determined grid (right), with contact contour.
(bottom) Zoomed in detail of boundary-determined grid (left) and operator-determined grid (right), with contact contour.}
\label{ex:obstacle:figC5}
\end{center}
\end{figure}

-  The free boundary-determined grid criteria specifies the finest grid resolution at nodes where both terms $F[u]$ and $u-g$ are close to zero.
\[
G^T[u]  \text{ determined by }  \min ( \abs{ \Delta u }, \abs{ u - g} )
\]
 The resulting grid provides the most refinement near the boundary of the contact set, as seen in Figure~\ref{ex:obstacle:figC5}.  

 - We also chose to refine the grid at nodes where the absolute value of operator $F^{obs}[u] = \min (-\Delta[u], u-g(x) )$ exceeds a threshold,
 \[
 G^F[u]  \text{ determined by } \abs{ F^{obs}[u] }. 
 \]
Notice in this case that the scaling of the two terms are different: the Laplacian scales like $1/h^2$ while the obstacle term has no scaling in $h$. 
  The resulting grid is very similar to the previous one, but with more refinement inside the contact set (which corresponds to capturing details of $g(x)$) , see Figure~\ref{ex:obstacle:figC5}.  
More Newton iterates were performed for the operator-determined grid, however these are often performed when the grid is still coarse resulting in overall slightly better performance.  

 The relative performance of the different refinement methods is given in Table~\ref{table:obstacle}.  

\begin{table}
\begin{center}
    \begin{tabular}{l |c|c|c|c }
    \multirow{2}{*}{Grid Type}&\multirow{2}{*}{Final nodes}&\multirow{2}{*}{Runtime}&\multicolumn{2}{| c }{Number of Newton Solves}\\ 
    \cline{4-5}
    &                          &                & Solves with $<5000$ nodes & total solves \\
    \hline
    Predetermined&22148&6.00s&21&59\\
    Boundary&15156&5.70s&25&67\\
    Operator&15967&5.54s&33&71
    \end{tabular}
    \caption{Comparison of adaptive grids for the obstacle problem.}
    \label{table:obstacle}	
\end{center}	
\end{table}

\subsection{Stefan free boundary problems}\label{ex:stefan}
We are now in the context of Example~\ref{ex:Stef}, where the Stefan problem is represented by the single operator $F^{Stef}$.

We take initial data corresponding to a function with three local maxima. 
To test the effectiveness of different grid adaptation strategies we solve the equation using:  \\
 -  a uniform coarse grid,\\ 
 -  a uniform fine grid,\\
 -  adapting the grid according to the size of the operator, $G^F[u] = \abs{F^{Stef[u]}}$,\\
 -  adapting the grid according to the size of both terms in the operator  $G^T[u] = \min (\abs{\Delta u}, \abs{\grad u}^2)$.
 
The reason for choosing the term adapted refinement $G^T$ comes from assuming that most of the accuracy of the solution comes from the accuracy of the free boundary: the term $G^T$ refines near the free boundary.   The operator-adapted refinement is intrinsic, it also scales correctly in terms of the grid, since both terms are order $1/h^2$.

\begin{figure}
\begin{center}
\includegraphics[keepaspectratio=false, width=\customheight, height=\customheight, angle=90]{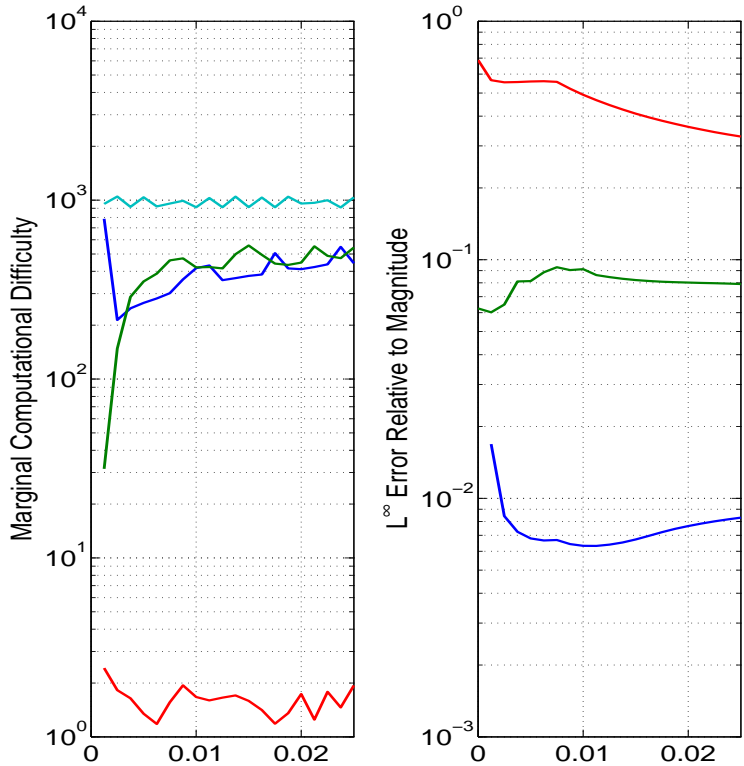}
\caption{Performance of the uniform-coarse grid (red), boundary-adaptive grid (green), operator-adaptive grid (blue) \& uniform-fine grid (light blue) when solving the sample Stefan equation.}
\label{ex:stefan:figA6}
\end{center}
\end{figure}

\begin{figure}
\begin{center}
\includegraphics[keepaspectratio=false, width=\customwidth, height=\customheight]{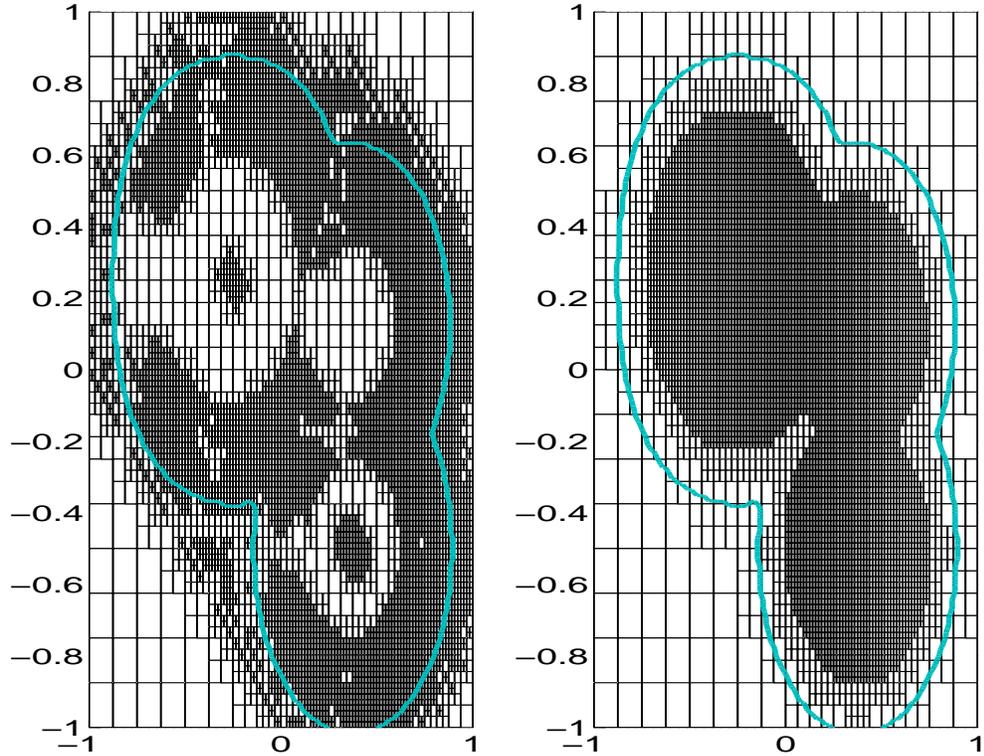}
\caption{Detail of the term-adapted (left) and operator-adapted (right) grids with the boundary contour of the uniform-fine grid solution of the sample Stefan equation at $t=0.005$.}
\label{ex:stefan:figC1}
\end{center}
\end{figure}

\begin{figure}
\begin{center}
\includegraphics[keepaspectratio=false, width=\customwidth, height=\customheight]{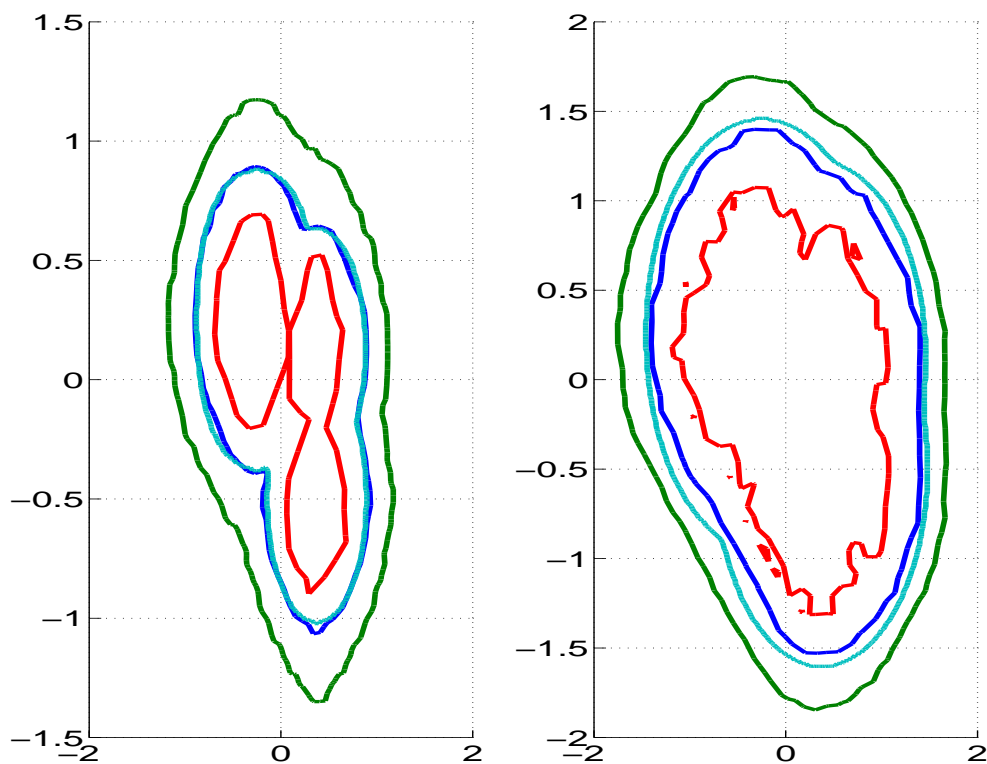}
\caption{Curves $\partial\{u=0\}$ of the sample Stefan problem at $t=0.005$ (left) and $t=0.025$ (right), due to solution with a uniform-coarse grid (red), term-adapted grid (green), operator-adapted grid (blue) \& uniform-fine grid (light blue).}
\label{ex:stefan:figB3}
\end{center}
\end{figure}

Figure \ref{ex:stefan:figA6} shows the observed computational complexity and outlines the $L^\infty$ proximity of the three methods to the fine-grid solution. By this metric, using $G^F$, the operator-adaptive grid, is superior to using $G^T$, the both terms adapted grid.  On closer inspection, Figure~\ref{ex:stefan:figC1} shows the $G^T$ grid is finest near $\partial\{u=0\}$, as desired, but that the solution is evolving slowly in this region. The operator $F^{Stef}$ is relatively small near the boundary.  By comparing the implied boundary-curves, Figure~\ref{ex:stefan:figB3} demonstrates clearly that the operator-adapted grid is a better strategy, because the accuracy of the location of the free boundary is significantly better in this case.

\section{Conclusions}

We  introduced a general framework for bringing adaptive grid and solvers to bear on a class of degenerate elliptic and parabolic Partial Differential Equations, which allows for the incorporation of free boundary problems, irregular and unbounded domains, along with adaptive grid refinement.  
We have demonstrated the significant improvement of solution accuracy and solution time, as compared to methods on regular grids.  

The adaptive grid overcomes the low accuracy of the finite difference method where curved boundaries are involved.  This includes the free boundaries which arise in the obstacle problem, or the Stefan problem.
By incorporating the boundary conditions etc into the operator, we can define a global residual, which included errors from the geometry as well as form the error.  Using this criteria we developed a grid refinement criteria which resulted in improvements over the other methods.

Other applications of the framework which are easily implemented include:
\begin{itemize}
\item First order equations, such as the eikonal equation in either bounded or unbounded domains.
\item Visibility problems, with refinement near essential small features of the obstructions 
\item Optimal or stochastic control problems with refinement at switching regions, so long as the second order operator can be discretized on the grid.
\item One phase free boundary problems, such as Hele-Shaw.
\end{itemize}

 
\bibliographystyle{alpha}
\bibliography{Adaptive}

\end{document}